**National and Kapodistrian University of Athens**

**Core Department**

**MSc in Advanced Control Systems and Robotics**

# «Supervisory Control of a Flexible Manufacturing Unit for the Production of Two Products»

Thesis

Kyriakos Giagiakos

# Abstract

In this diploma thesis, the mathematical model of a multi-product manufacturing unit will be presented. The unit consists of a set of three conveyors, a robot, a lathe, a milling machine, an assembly machine, and a painting machine. Finally, the connection of the above elements is carried out via one slot buffer. The products can be divided into two distinct sets according to the route they will follow through the machines of the unit. The mathematical models of the individual subsystems of the plant using finite deterministic automata will be presented. The two different routes that the products can follow will be presented. The desired product flow will be presented in the form of desired regular languages. The properties of the desired languages with respect the overall automaton of the system will be investigated. A supervisory architecture will be designed based on the desired regular languages.

# CHAPTER 1: INTRODUCTION

Flexible manufacturing systems are quite frequently encountered in modern industry. Indicative examples of such systems are presented in works [1] and [2]. The control of these systems requires the design of control architectures that can easily adapt to possible changes in both the processing sequence of the products and the system's layout itself, see [1]-[3]. In work [3], a flexible manufacturing system is presented, which has two input lines and can produce two different products. Extending this problem to multiple input lines for the production of two different products using discrete event supervisors is currently under investigation.

In this thesis, the mathematical model of a multi-product processing unit will be presented. The unit consists of three conveyors, a robot, a lathe, a milling machine, an assembly machine, and a painting machine. The interconnection of the above components is implemented through temporary storage buffers. The products can be divided into two distinct sets, depending on the path they follow through the machines of the unit. The mathematical models of the system's individual subsystems will be presented using finite deterministic automata. The two different processing paths that the products can follow will be analyzed. The desired product flows will be expressed in terms of regular languages. The properties of these desired languages will be examined in relation to the complete automaton of the system. A supervisory control architecture will be designed.

A detailed description of the structure of the thesis is presented below:

**Chapter 2** presents the mathematical model of the processing unit and its subsystems using finite deterministic automata.

**Chapter 3** describes the desired behavior of the system in the form of rules and regular languages. It also proves the controllability of the desired languages with respect to the overall system.

**Chapter 4** presents the design of the supervisory control architecture.

The thesis concludes with a summary of the main findings.

# CHAPTER 2: Modeling and Simulation of the Industrial Process System

## 2.1 Architecture of the system

Figure 1 presents the structure of the processing system for two products. The two products enter the processing system through two different input points (conveyor belts)."

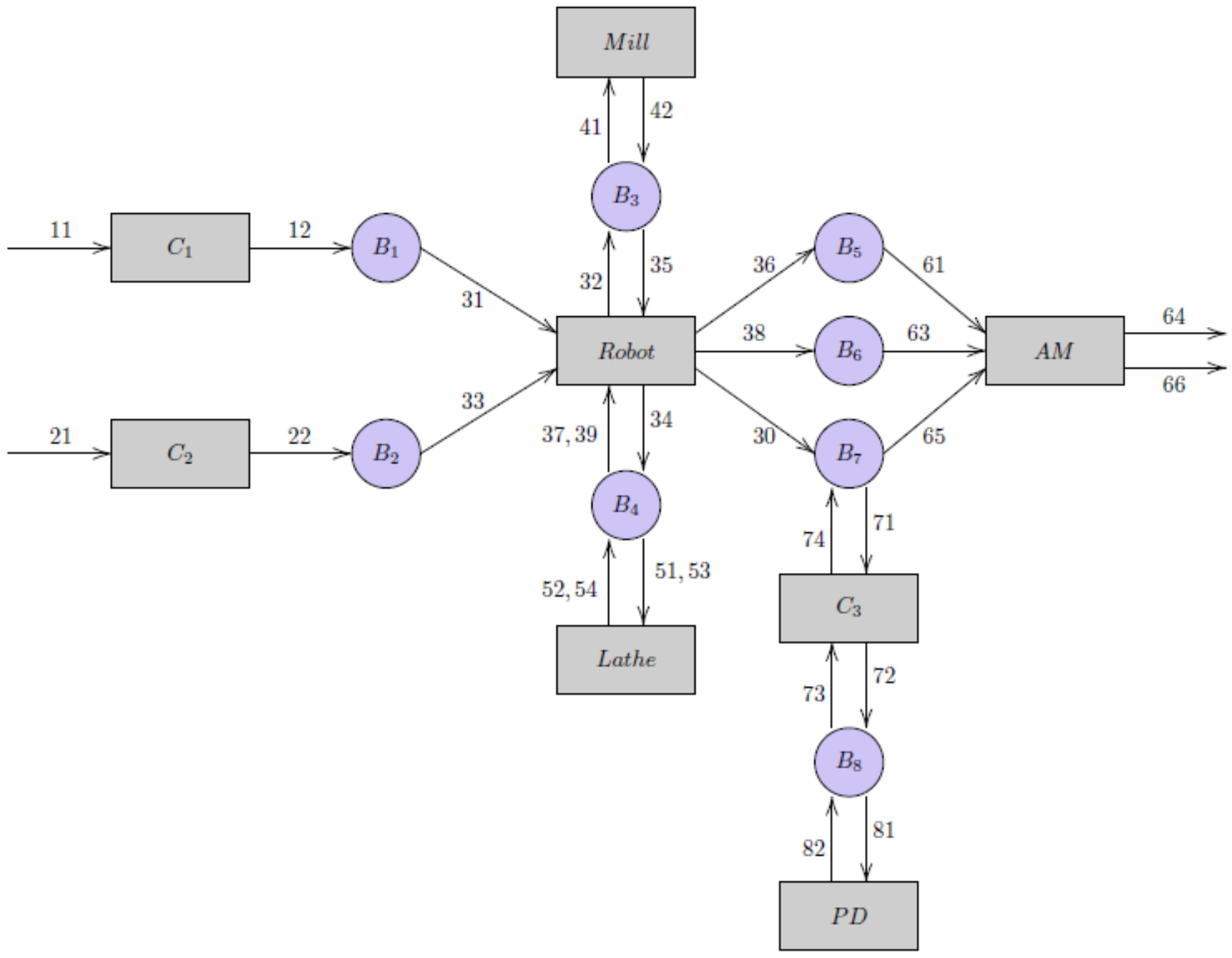

Figure 1: Architecture of the system *[1]*

The processing system consists of eight subsystems: three conveyor belts, a milling machine, a lathe, a transfer robot, a painting device, and an assembly machine, as well as eight temporary storage positions (buffers). All conveyor belts are denoted by the letter $\mathbf{C}_i$, $i = \{1, 2, 3\}$. All buffers are denoted with the latter $\mathbf{B}_j$, $j = \{1, 2, 3, 4, 5, 6, 7, 8\}$. The milling machine is denoted by the letter $\mathbf{M}$. The lathe is denoted by the letter $\mathbf{L}$. The robotic maniplulator is denoted by the letter $\mathbf{R}$. The painting machine is denoted by the letter $\mathbf{P}$. Finally, the assembly machine is denoted by the letter $\mathbf{A}$. The flow of the products is as follows:

The products enter the system via the conveyor belts $\mathbf{C}_1$ and $\mathbf{C}_2$ into the temporary buffer positions $\mathbf{B}_1$ and $\mathbf{B}_2$, respectively. The robotic arm $\mathbf{R}$ transfers and receives products to and from the temporary buffers $\mathbf{B}_1$ to $\mathbf{B}_7$. Buffer $\mathbf{B}_3$ connects the robotic arm with the milling machine. Buffer $\mathbf{B}_4$ connects the robotic arm with the lathe. Buffers $\mathbf{B}_5$ and $\mathbf{B}_6$ connect the robotic arm with the assembly machine. Buffer $\mathbf{B}_7$ connects the robotic arm with conveyor belt $\mathbf{C}_3$ and with the assembly machine. Buffer $\mathbf{B}_8$ connects conveyor belt $\mathbf{C}_3$ with the painting machine. In the following subsections, the mathematical models of all processing system components will be presented.

## 2.2 Model of the three conveyor belts

The conveyor belt performs discrete movements. The automaton of conveyor belt $C_i$, where $i = \{1,2,3\}$, is described by the 6-tuple ([4]-[16]).

$$\mathbf{G}_{C,i} = (\mathbb{Q}_{C,i}, \mathbb{E}_{C,i}, f_{C,i}, \mathbb{H}_{C,i}, x_{C,i,0}, \mathbb{Q}_{C,i,m})$$

The set $\mathbb{Q}_{C,i} = \{q_{C,i,1}, q_{C,i,2}\}$ is the set of states of the conveyor belt. In state $q_{C,i,1}$, the conveyor belt is empty. In state $q_{\mathrm{C}i,2}$, the conveyor belt is transporting a product.

The initial state of the conveyor belt is $x_{C,i,0} = q_{C,i,1}$, which is the state where the conveyor belt is stationary. The set of marked states of the conveyor belt is $\mathbb{Q}_{C,i,m} = \{q_{C,i,1}\}$.

The alphabet of the conveyor belt is $\mathbb{E}_{C,i} = \{e_{C,i,1}, e_{C,i,2}\}$. The event $e_{C,i,1}$ corresponds to the signal that a product has been loaded. The event $e_{\mathrm{C}i,2}$ corresponds to the command to move the conveyor belt. The set of controllable events of the system is $\mathbb{E}_{C,i,c} = \{e_{C,i,2}\}$, and the set of uncontrollable events of the system is $\mathbb{E}_{C,i,uc} = \{e_{C,i,2}\}$.

The transitions caused by each event of the conveyor belt are

$$f_{C,i}(q_{C,i,1}, e_{C,i,1}) = q_{C,i,2} \text{ and } f_{C,i}(q_{C,i,2}, e_{C,i,2}) = q_{C,i,1}$$

The sets of active events for each state are

$$\mathbb{H}_{C,i}(q_{C,i,1}) = \{e_{C,i,1}\} \text{ and } \mathbb{H}_{C,i}(q_{C,i,2}) = \{e_{C,i,2}\}$$

The language produced by the conveyor belt automaton is

$$\mathbb{L}(\mathbf{G}_{C,i}) = \overline{(e_{C,i,1} e_{C,i,2})^*}$$

Finally, the language marked by the conveyor belt automaton is

$$\mathbb{L}_m(\mathbf{G}_{C,i}) = (e_{C,i,1}e_{C,i,2})^*$$

The state diagram of the conveyor belt is presented in Figure 2

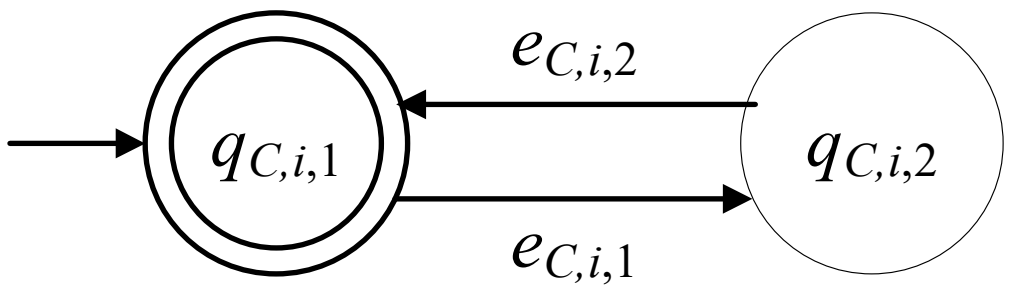

.

Figure 2: State diagram of the *i*-th conveyor belt

## 2.3 Model of the robotic manipulator

Next, the automaton of the robotic arm $\mathbf{R}$ will be described using finite deterministic automata. In [17], the dynamic description of a robotic arm is presented, along with the commands for moving the arm to and from predefined positions. Therefore, the automaton $\mathbf{R}$ is described by the 6-tuple ([4]-[16]).

$$\mathbf{G}_R = (\mathbb{Q}_R, \mathbb{E}_R, f_R, \mathbb{H}_R, x_{R,0}, \mathbb{Q}_{R,m})$$

The set of states of the robotic arm is $\mathbb{Q}_R = \{q_{R,1}, q_{R,2}\}$. State $q_{R,1}$ is the state where the robotic arm is empty. State $q_{R,2}$ is the state where the robotic arm is full.

The initial state of the system is $x_{R,0} = q_{R,1}$, which is the state where the robotic arm is empty. The set of marked states is $\mathbb{Q}_{R,m} = \{q_{R,1}\}$.

The alphabet is

$$\mathbb{E}_R = \{e_{R,1}, e_{R,2}, e_{R,3}, e_{R,4}, e_{R,5}, e_{R,6}, e_{R,7}, e_{R,8}, e_{R,9}, e_{R,10}, e_{R,11}, e_{R,12}, e_{R,13}, e_{R,14}\}$$

The event $e_{R,\lambda}$ corresponds to the command to load a product onto the robotic arm from the temporary memory position $\mathbf{B}_\lambda$, where $\lambda \in \{1,2,3,4,5,6,7\}$. The event $e_{R,\lambda+7}$ corresponds to the command to place a product in the temporary memory position $\mathbf{B}_\lambda$ by the robotic arm, where $\lambda \in \{1,2,3,4,5,6,7\}$. All events of the robotic arm $\mathbf{R}$ are controllable, thus it holds that $\mathbb{E}_{R,c} = \mathbb{E}_R$ and $\mathbb{E}_{R,uc} = \varnothing$.

The transitions caused by each event of the robotic arm are

$$f_R(q_{R,1}, e_{R,\lambda}) = q_{\mathrm{R1},2}\,,\ \lambda \in \{1,2,3,4,5,6,7\}$$

$$f_R(q_{R,2}, e_{R,\lambda+7}) = q_{\mathrm{R1,1}}\,,\ \lambda \in \{1,2,3,4,5,6,7\}$$

The sets of active events for each state are

$$\mathbb{H}_R(q_{R,1}) = \{e_{R,1}, e_{R,2}, e_{R,3}, e_{R,4}, e_{R,5}, e_{R,6}, e_{R,7}\} \text{ and}$$
$$\mathbb{H}_R(q_{R,2}) = \{e_{R,8}, e_{R,9}, e_{R,10}, e_{R,11}, e_{R,12}, e_{R,13}, e_{R,14}\}$$

The language produced by the robotic arm automaton is

$$\mathbb{L}(\mathbf{G}_R) = \overline{\Big((e_{R,1} + e_{R,2} + e_{R,3} + e_{R,4} + e_{R,5} + e_{R,6} + e_{R,7}) (e_{R,8} + e_{R,9} + e_{R,10} + e_{R,11} + e_{R,12} + e_{R,13} + e_{R,14})\Big)^*}$$

Finally, the language marked by the robotic arm automaton is

$$\mathbb{L}(\mathbf{G}_R) = \Big((e_{R,1} + e_{R,2} + e_{R,3} + e_{R,4} + e_{R,5} + e_{R,6} + e_{R,7}) (e_{R,8} + e_{R,9} + e_{R,10} + e_{R,11} + e_{R,12} + e_{R,13} + e_{R,14})\Big)^*$$

The state diagram of the robotic arm is presented in Figure 3.

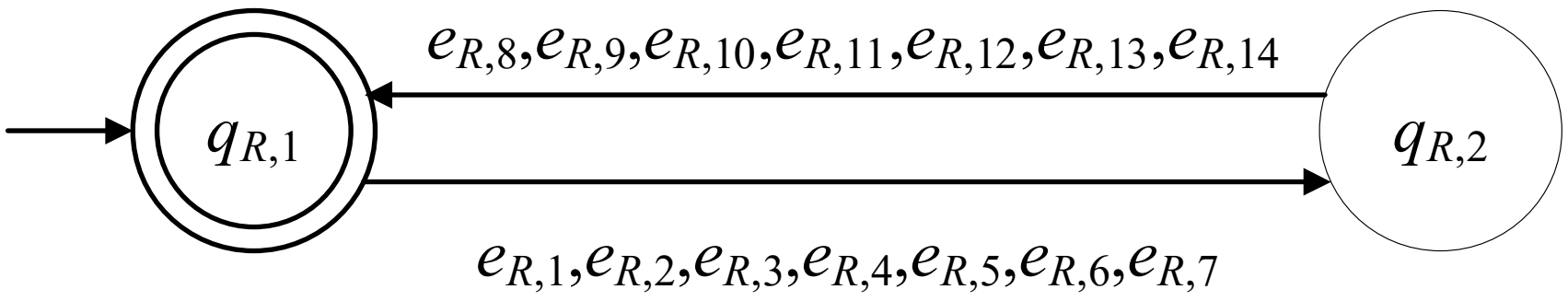

Figure 3: Automaton of the robotic arm

## 2.4 Model of the lathe

The automaton of the lathe is

$$\mathbf{G}_L = (\mathbb{Q}_L, \mathbb{E}_L, f_L, \mathbb{H}_L, x_{L,0}, \mathbb{Q}_{L,m})$$

The set of states of the lathe is $\mathbb{Q}_L = \{q_{L,1}, q_{L,2}\}$. State $q_{L,1}$ is the state where the lathe is turned off. State $q_{L,2}$ is the state where the lathe is turned on.

The initial state of the lathe is $x_{L,0} = q_{L,1}$, which is the state where the lathe is turned off. The set of marked states of the lathe is $\mathbb{Q}_{L,m} = \{q_{L,1}\}$.

The alphabet of the lathe is $\mathbb{E}_L = \{e_{L,1}, e_{L,2}, e_{L,3}, e_{L,4}\}$. Event $e_{L,1}$ is the command to load a product of the first category and start the machining process. Event $e_{L,2}$ is the command to load a product of the second category and start the machining process. Event $e_{L,3}$ is the command to complete the machining of the product of the first category. Event $e_{L,4}$ is the command to complete the machining of the product of the second category. All events of the lathe are controllable, thus it holds that $\mathbb{E}_{L,c} = \mathbb{E}_L$ and $\mathbb{E}_{L,uc} = \varnothing$.

The transitions caused by each event of the lathe are.

$$f_L(q_{L,1}, e_{L,1}) = q_{L,2},\ f_L(q_{L,1}, e_{L,2}) = q_{L,2}$$

$$f_L(q_{L,2}, e_{L,3}) = q_{L,1} \text{ and } f_L(q_{L,2}, e_{L,4}) = q_{L,1}$$

The sets of active events for each state of the lathe are

$$\mathbb{H}_L(q_{L,1}) = \{e_{L,1}, e_{L,2}\} \text{ and } \mathbb{H}_L(q_{L,2}) = \{e_{L,3}, e_{L,4}\}$$

The language produced by the lathe automaton is

$$\mathbb{L}(\mathbf{G}_L) = \overline{\left((e_{L,1} + e_{L,2})(e_{L,3} + e_{L,4})\right)^*}$$

Finally, the language marked by the lathe automaton is

$$\mathbb{L}_m(\mathbf{G}_L) = \left((e_{L,1} + e_{L,2})(e_{L,3} + e_{L,4})\right)^*$$

The state diagram of the lathe is presented in Figure 4.

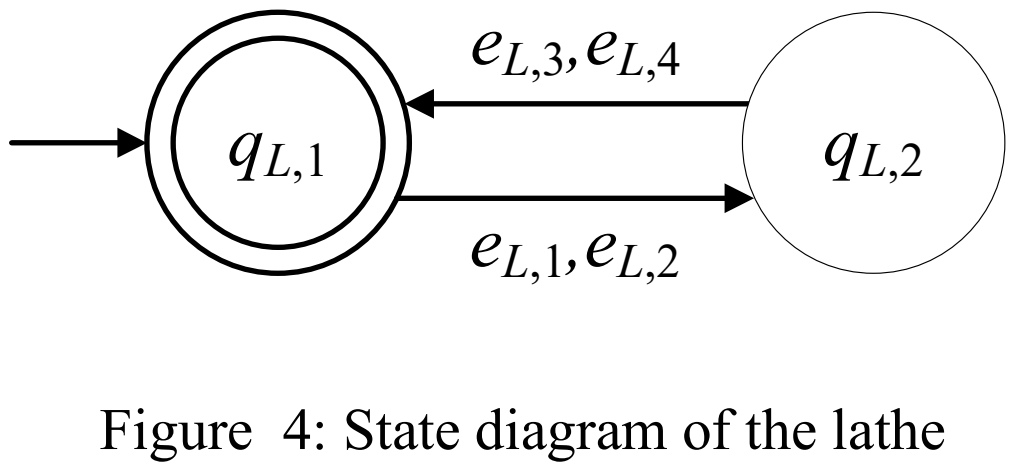

Figure 4: State diagram of the lathe

## 2.5 Model of the milling machine

The automaton of the milling machine is

$$\mathbf{G}_M = (\mathbb{Q}_M, \mathbb{E}_M, f_M, \mathbb{H}_M, x_{M,0}, \mathbb{Q}_{M,m})$$

The set of states of the milling machine is $\mathbb{Q}_M = \{q_{M,1}, q_{M,2}\}$. State $q_{M,1}$ is the state where the milling machine is turned off. State $q_{M,2}$ is the state where the milling machine is turned on.

The initial state of the milling machine is $x_{M,0} = q_{M,1}$, which is the state where the milling machine is turned off. The set of marked states of the milling machine is $\mathbb{Q}_{M,m} = \{q_{M,1}\}$.

The alphabet of the milling machine is $\mathbb{E}_M = \{e_{M,1}, e_{M,2}\}$. Event $e_{M,1}$ is the command to load a product and start the machining process. Event $e_{M,2}$ is the command to complete the machining of the product. All events of the milling machine are controllable, thus it holds that $\mathbb{E}_{M,c} = \mathbb{E}_M$ and $\mathbb{E}_{M,uc} = \varnothing$.

The transitions caused by each event of the milling machine are

$$f_M(q_{M,1}, e_{M,1}) = q_{M,2} \text{ and } f_M(q_{M,2}, e_{M,2}) = q_{M,1}$$

The sets of active events for each state of the milling machine are

$$\mathbb{H}_M(q_{M,1}) = \{e_{M,1}\} \text{ and } \mathbb{H}_M(q_{M,2}) = \{e_{M,2}\}$$

The language produced by the milling machine automaton is

$$\mathbb{L}(\mathbf{G}_M) = \overline{(e_{M,1}e_{M,2})^*}$$

Finally, the language marked by the milling machine automaton is

$$\mathbb{L}_m(\mathbf{G}_M) = (e_{M,1}e_{M,2})^*$$

The state diagram of the milling machine is presented in Figure 5.

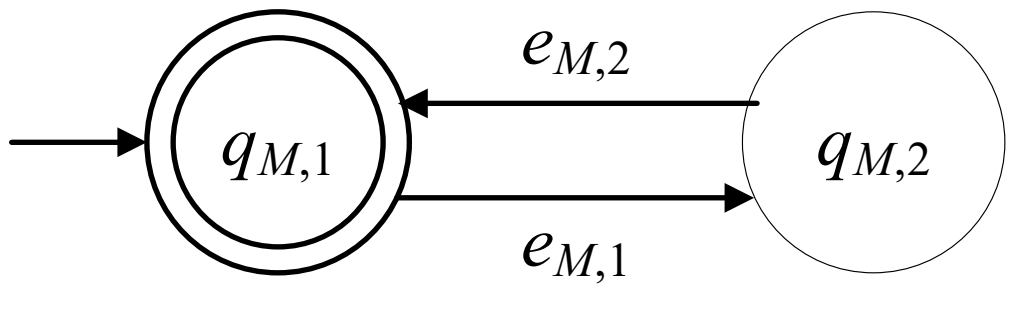

Figure 5: Automaton of the milling machine

## 2.6 Model of the painting machine

The automaton of the painting machine is

$$\mathbf{G}_P = (\mathbb{Q}_P, \mathbb{E}_P, f_P, \mathbb{H}_P, x_{P,0}, \mathbb{Q}_{P,m}).$$

The set of states of the painting machine is $\mathbb{Q}_P = \{q_{P,1}, q_{P,2}\}$. State $q_{P,1}$ is the state where the painting machine is turned off. State $q_{P,2}$ is the state where the painting machine is turned on.

The initial state of the painting machine is $x_{P,0} = q_{P,1}$, which is the state where the painting machine is turned off. The set of marked states of the painting machine is $\mathbb{Q}_{P,m} = \{q_{P,1}\}$.

The alphabet of the painting machine is $\mathbb{E}_P = \{e_{P,1}, e_{P,2}\}$. Event $e_{P,1}$ is the command to load a product and start the painting process. Event $e_{P,2}$ is the command to complete the painting of the product. All events of the painting machine are controllable, thus it holds that $\mathbb{E}_{P,c} = \mathbb{E}_P$ and $\mathbb{E}_{P,uc} = \varnothing$.

The transitions caused by each event of the painting machine are.

$$f_P(q_{P,1}, e_{P,1}) = q_{P,2} \text{ and } f_P(q_{P,2}, e_{P,2}) = q_{P,1}$$

The sets of active events for each state of the painting machine are

$$\mathbb{H}_P(q_{P,1}) = \{e_{P,1}\} \text{ and } \mathbb{H}_P(q_{P,2}) = \{e_{P,2}\}$$

The produced language by the painting machine automaton is

$$\mathbb{L}(\mathbf{G}_P) = \overline{(e_{P,1}e_{P,2})^*}$$

Finally, the marked language by the painting machine automaton is

$$\mathbb{L}_m(\mathbf{G}_P) = (e_{P,1}e_{P,2})^*$$

The state diagram of the painting machine is presented in Figure 6.

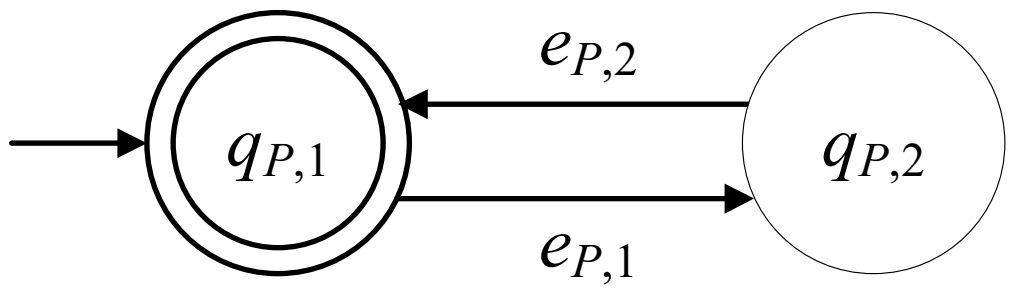

Figure 6: Automaton of the painting machine

## 2.7 Model of the assembly machine

The automaton of the assembly machine is

$$\mathbf{G}_A = (\mathbb{Q}_A, \mathbb{E}_A, f_A, \mathbb{H}_A, x_{A,0}, \mathbb{Q}_{A,m})$$

The set of states of the milling machine is $\mathbb{Q}_A = \{q_{A,1}, q_{A,2}, q_{A,3}\}$. State $q_{A,1}$ is the state where the assembly machine is turned off. State $q_{A,2}$ is the state where the assembly machine is in stand by mode. State $q_{A,3}$ is the state where the assembly machine is turned on.

The initial state of the assembly machine is $x_{A,0} = q_{A,1}$, which is the state where the assembly machine is turned off. The set of marked states of the assembly machine is $\mathbb{Q}_{A,m} = \{q_{A,1}\}$.

The alphabet of the assembly machine is $\mathbb{E}_A = \{e_{A,1}, e_{A,2}, e_{A,3}, e_{A,4}, e_{A,5}, e_{A,6}\}$. Event $e_{A,1}$ is the command to start the assembly machine. Event $e_{A,2}$ is the command to assemble a product from the temporary memory position $\mathbf{B}_5$. Event $e_{A,3}$ is the command to assemble a product from the temporary memory position $\mathbf{B}_6$. Event $e_{A,4}$ is the command to assemble a product from the temporary memory position $\mathbf{B}_7$. Event $e_{A,5}$ is the command that the assembly of the first category products has been completed. Event $e_{A,6}$ is the command that the assembly of the second category products has been completed. Events $e_{A,1}$, $e_{A,2}$, $e_{A,3}$ and $e_{A,4}$ are controllable events. Therefore, the set of controllable events is $\mathbb{E}_{A,c} = \{e_{A,1}, e_{A,2}, e_{A,3}, e_{A,4}\}$. The set of uncontrollable events is $\mathbb{E}_{A,uc} = \{e_{A,5}, e_{A,6}\}$.

The transitions caused by each event of the assembly machine are

$$f_A(q_{A,1}, e_{A,1}) = q_{A,2},\ f_A(q_{A,2}, e_{A,2}) = q_{A,3},$$

$$f_A(q_{A,2}, e_{A,3}) = q_{A,3},\ f_A(q_{A,2}, e_{A,4}) = q_{A,3}$$

$$f_A(q_{A,3}, e_{A,5}) = q_{A,1} \text{ and } f_A(q_{A,3}, e_{A,6}) = q_{A,1}$$

The sets of active events for each state of the assembly machine are

$$\mathbb{H}_A(q_{A,1}) = \{e_{A,1}\},\ \mathbb{H}_A(q_{A,2}) = \{e_{A,2}, e_{A,3}, e_{A,4}\}$$
$$\text{and } \mathbb{H}_A(q_{A,3}) = \{e_{A,5}, e_{A,6}\}$$

The produced language by the assembly machine automaton is

$$\mathbb{L}(\mathbf{G}_A) = \overline{\left(e_{A,1}(e_{A,2} + e_{A,3} + e_{A,4})(e_{A,5} + e_{A,6})\right)^*}$$

Finally, the marked language by the assembly machine automaton is

$$\mathbb{L}_m(\mathbf{G}_A) = \left(e_{A,1}(e_{A,2} + e_{A,3} + e_{A,4})(e_{A,5} + e_{A,6})\right)^*$$

The state diagram of the assembly machine is presented in Figure 7.

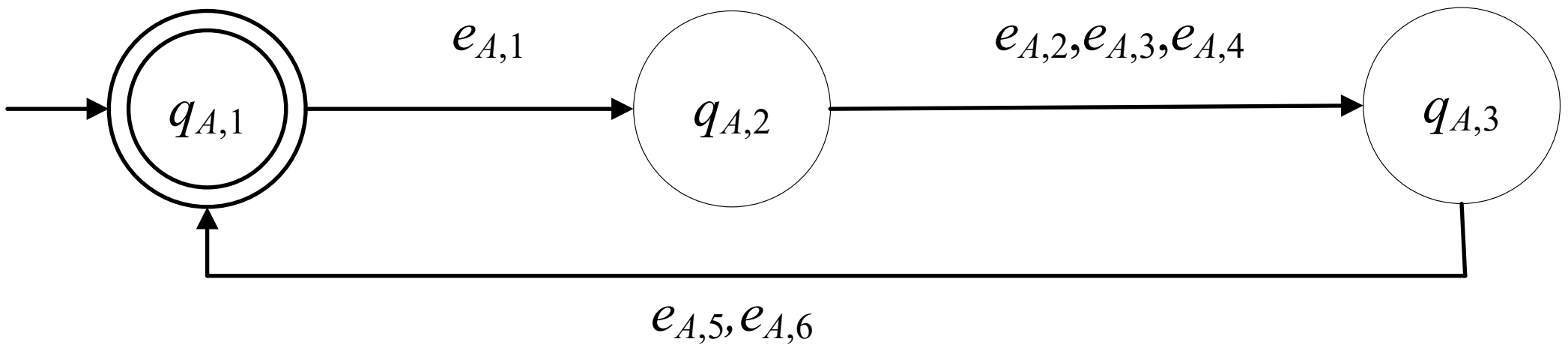

Figure 7: Automaton of the assembly machine

## 2.8 The total model

The overall automaton is described by the parallel composition of all the above automata (a total of 8 automata) (more on parallel composition can be found in references [18] and [19]). The resulting overall automaton is

$$\mathbf{G} = \mathbf{G}_{C,1} \,||\, \mathbf{G}_{C,2} \,||\, \mathbf{G}_{C,3} \,||\, \mathbf{G}_R \,||\, \mathbf{G}_L \,||\, \mathbf{G}_M \,||\, \mathbf{G}_P \,||\, \mathbf{G}_A$$

It holds that

$$\mathbb{Q} = \mathbb{Q}_{C,1} \times \mathbb{Q}_{C,2} \times \mathbb{Q}_{C,3} \times \times \mathbb{Q}_{R} \times \mathbb{Q}_{L} \times \mathbb{Q}_{M} \times \mathbb{Q}_{P} \times \mathbb{Q}_{A}$$

That is, the states will be described by the 8-tuple

$$(q_{C,1}, q_{C,2}, q_{C,3}, q_{R}, q_{L}, q_{M}, q_{P}, q_{A}),$$

The events of the system are

$$\mathbb{E} = \mathbb{E}_{C,1} \cup \mathbb{E}_{C,2} \cup \mathbb{E}_{C,3} \cup \mathbb{E}_{R} \cup \mathbb{E}_{L} \cup \mathbb{E}_{M} \cup \mathbb{E}_{P} \cup \mathbb{E}_{A}$$

The controllable alphabet of the system is

$$\mathbb{E}_{c} = \mathbb{E}_{C,1,c} \cup \mathbb{E}_{C,2,c} \cup \mathbb{E}_{C,3,c} \cup \mathbb{E}_{R,c} \cup \mathbb{E}_{L,c} \cup \mathbb{E}_{M,c} \cup \mathbb{E}_{P,c} \cup \mathbb{E}_{A,c}$$

The uncontrollable alphabet of the system is

$$\begin{aligned} \mathbb{E}_{c} &= \mathbb{E}_{C,1,uc} \cup \mathbb{E}_{C,2,uc} \cup \mathbb{E}_{C,3,uc} \cup \mathbb{E}_{R,uc} \cup \mathbb{E}_{L,uc} \cup \mathbb{E}_{M,uc} \cup \mathbb{E}_{P,uc} \cup \mathbb{E}_{A,uc} \\ &= \{e_{C,1,2}, e_{C,2,2}, e_{C,3,3}, e_{A,5}, e_{A,6}\} \end{aligned}$$

The initial state is

$$x_{0} = (q_{C,1,1}, q_{C,2,1}, q_{C,3,1}, q_{R,1}, q_{L,1}, q_{M,1}, q_{P,1}, q_{A,1})$$

The set of marked states is

$$\mathbb{Q}_{m} = \left\{(q_{C,1,1}, q_{C,2,1}, q_{C,3,1}, q_{R,1}, q_{L,1}, q_{M,1}, q_{P,1}, q_{A,1})\right\}$$

# CHAPTER 3: DESIRED BEHAVIOR OF THE INDUSTRIAL SYSTEM AND REGULAR LANGUAGES

## 3.1 Desired behavior in the form of rules

The desired behavior of an industrial system can be described and designed using static or dynamic controllers in continuous or discrete time ([20]), or in the form of regular languages ([18]–[19]). The description of the desired behavior is as follows:

**Manufacturing of the products**:

- Category 1 products enter through conveyor belt $\mathbf{C}_1$ and are stored in temporary buffer $\mathbf{B}_1$. Then, they are sequentially processed by the milling machine, the lathe, the painting machine (for certain products), and the assembly machine. Therefore, the product flow is

$$\mathbf{C}_1 \rightarrow \mathbf{B}_1 \rightarrow \mathbf{R} \rightarrow \mathbf{B}_3 \rightarrow \mathbf{M} \rightarrow \mathbf{B}_3 \rightarrow \mathbf{R} \rightarrow \mathbf{B}_4 \rightarrow \mathbf{L} \rightarrow \mathbf{B}_4 \rightarrow \mathbf{R} \rightarrow \mathbf{B}_6 \rightarrow \mathbf{A}$$
$$\downarrow$$
$$\mathbf{A} \leftarrow \mathbf{B}_7 \leftarrow \mathbf{C}_3 \leftarrow \mathbf{B}_8 \leftarrow \mathbf{P} \leftarrow \mathbf{B}_8 \leftarrow \mathbf{C}_3 \leftarrow \mathbf{B}_7$$

- Category 2 products enter through conveyor belt $\mathbf{C}_2$ and are stored in temporary buffer $\mathbf{B}_2$. Then, they are sequentially processed by the lathe, the painting machine (for certain products), and the assembly machine. Therefore, the product flow is

$$\mathbf{C}_2 \rightarrow \mathbf{B}_2 \rightarrow \mathbf{R} \rightarrow \mathbf{B}_4 \rightarrow \mathbf{L} \rightarrow \mathbf{B}_4 \rightarrow \mathbf{R} \rightarrow \mathbf{B}_5 \rightarrow \mathbf{A}$$
$$\downarrow$$
$$\mathbf{A} \leftarrow \mathbf{B}_7 \leftarrow \mathbf{C}_3 \leftarrow \mathbf{B}_8 \leftarrow \mathbf{P} \leftarrow \mathbf{B}_8 \leftarrow \mathbf{C}_3 \leftarrow \mathbf{B}_7$$

## 3.2 Desired behavior in the form of regular languages

### *3.2.1 Regular language for the first category*

The regular language that describes the desired behavior for the first product is of the form:

$$\mathbb{K}_{D,1} = \overline{\left(e_{C,1,1}e_{R,1}e_{R,10}e_{M,1}e_{R,3}e_{R,11}e_{L,1}e_{R,4}(e_{R,13} + e_{R,14}e_{C,3,1}e_{P,1}e_{C,3,1})e_{A,1}\right)^{*}}$$

The controllability of language $\mathbb{K}_{D,1}$ with respect to the overall automaton $\mathbf{G}$ is ensured by the fact that only controllable events participate in the language.

### *3.2.2 Regular language for the second category*

The regular language that describes the desired behavior for the second product is of the form:

$$\mathbb{K}_{D,2} = \overline{\left(e_{C,2,1}e_{R,2}e_{R,11}e_{L,2}e_{R,4}(e_{R,12} + e_{R,14}e_{C,3,1}e_{P,1}e_{C,3,1})e_{A,1}\right)^{*}}$$

The controllability of language $\mathbb{K}_{D,2}$ with respect to the overall automaton $\mathbf{G}$ is ensured by the fact that only controllable events participate in the language.

# CHAPTER 4: SUPERVISORS

## 4.1 Supervisor for the first category

For the desired language $\mathbb{K}_{D,1}$, the supervisor's automaton will be of the form

$$\mathbf{S}_1 = (\mathbb{Q}_{S,1}, \mathbb{E}_{S,1}, f_{S,1}, \mathbb{H}_{S,1}, x_{S,1,0}, \mathbb{Q}_{S,1,m})$$

The set of states of the supervisor is

$$\mathbb{Q}_{S,1} = \bigcup_{i=1}^{14} \{q_{S,1,i}\}$$

The alphabet of the supervisor is

$$\mathbb{E}_{S,1} = \{e_{C,1,1}, e_{C,3,1}, e_{R,1}, e_{R,3}, e_{R,4}, e_{R,11}, e_{R,10}, e_{R,13}, e_{R,14}, e_{M,1}, e_{L,1}, e_{P,1}, e_{A,1}\}$$

The transition functions are

$$f_{S,1}(q_{S,1,1}, e_{C,1,1}) = q_{S,1,2},\ f_{S,1}(q_{S,1,2}, e_{R,1}) = q_{S,1,3},$$

$$f_{S,1}(q_{S,1,3}, e_{R,10}) = q_{S,1,4},\ f_{S,1}(q_{S,1,4}, e_{M,1}) = q_{S,1,5},$$

$$f_{S,1}(q_{S,1,5}, e_{R,3}) = q_{S,1,6},\ f_{S,1}(q_{S,1,6}, e_{R,11}) = q_{S,1,7},$$

$$f_{S,1}(q_{S,1,7}, e_{L,1}) = q_{S,1,8},\ f_{S,1}(q_{S,1,8}, e_{R,4}) = q_{S,1,9},$$

$$f_{S,1}(q_{S,1,9}, e_{R,13}) = q_{S,1,10},\ f_{S,1}(q_{S,1,9}, e_{R,14}) = q_{S,1,11},$$

$$f_{S,1}(q_{S,1,10}, e_{A,1}) = q_{S,1,1},\ f_{S,1}(q_{S,1,11}, e_{C,3,1}) = q_{S,1,12},$$

$$f_{S,1}(q_{S,1,12}, e_{P,1}) = q_{S,1,13},\ f_{S,1}(q_{S,1,13}, e_{C,3,1}) = q_{S,1,14}$$

and $f_{S,1}(q_{S,1,14}, e_{A,1}) = q_{S,1,1}$.

The active events sets are

$$\mathbb{H}_{S,1}(q_{S,1,1}) = \{e_{C,1,1}\} \ \mathbb{H}_{S,1}(q_{S,1,2}) = \{e_{R,1}\},$$

$$\mathbb{H}_{S,1}(q_{S,1,3}) = \{e_{R,10}\},\ \mathbb{H}_{S,1}(q_{S,1,4}) = \{e_{M,1}\},$$

$$\mathbb{H}_{S,1}(q_{S,1,5}) = \{e_{R,3}\},\ \mathbb{H}_{S,1}(q_{S,1,6}) = \{e_{R,11}\},$$

$$\mathbb{H}_{S,1}(q_{S,1,7}) = \{e_{L,1}\},\ \mathbb{H}_{S,1}(q_{S,1,8}) = \{e_{R,4}\},$$

$$\mathbb{H}_{S,1}(q_{S,1,9}) = \{e_{R,13}, e_{R,14}\},\ \mathbb{H}_{S,1}(q_{S,1,10}) = \{e_{A,1}\},$$

$$\mathbb{H}_{S,1}(q_{S,1,11}) = \{e_{C,3,1}\},\ \mathbb{H}_{S,1}(q_{S,1,12}) = \{e_{P,1}\},$$

$$\mathbb{H}_{S,1}(q_{S,1,13}) = \{e_{C,3,1}\} \text{ and } \mathbb{H}_{S,1}(q_{S,1,14}) = \{e_{A,1}\}.$$

The initial state is $x_{S,1,0} = q_{S,1,1}$.

The set of marked states is $\mathbb{Q}_{S,1,m} = \mathbb{Q}_{S,1}$.

The language generated by the automaton is

$$\mathbb{L}(\mathbf{S}_1) = \overline{\left(e_{C,1,1}e_{R,1}e_{R,10}e_{M,1}e_{R,3}e_{R,11}e_{L,1}e_{R,4}(e_{R,13} + e_{R,14}e_{C,3,1}e_{P,1}e_{C,3,1})e_{A,1}\right)^*}$$

Finally, the marked language of the automaton is

$$\mathbb{L}_m(\mathbf{S}_1) = \overline{\left(e_{C,1,1}e_{R,1}e_{R,10}e_{M,1}e_{R,3}e_{R,11}e_{L,1}e_{R,4}(e_{R,13} + e_{R,14}e_{C,3,1}e_{P,1}e_{C,3,1})e_{A,1}\right)^*}$$

The supervisor’s automaton is shown in Figure 8.

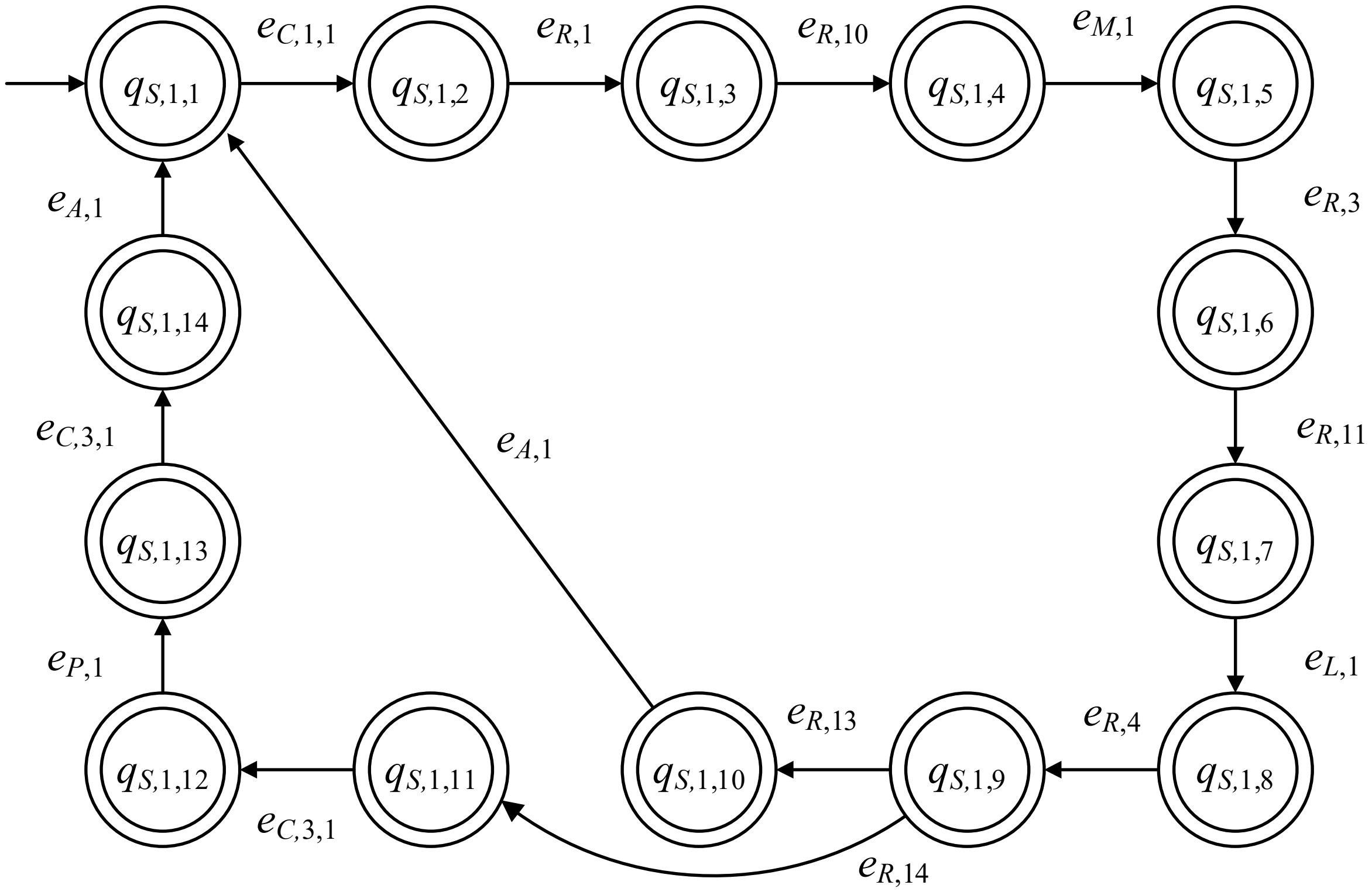

Figure 8: Automation of $\mathbf{S_1}$

## 4.2 Supervisor for the second category

For the desired language $\mathbb{K}_{D,2}$, the supervisor's automaton will be of the form

$$\mathbf{S}_2 = (\mathbb{Q}_{S,2}, \mathbb{E}_{S,2}, f_{S,2}, \mathbb{H}_{S,2}, x_{S,2,0}, \mathbb{Q}_{S,2,m})$$

The set of states of the supervisor is

$$\mathbb{Q}_{S,2} = \bigcup_{i=1}^{11} \{q_{S,2,i}\}$$

The alphabet of the supervisor is

$$\mathbb{E}_{S,2} = \{e_{C,2,1}, e_{C,3,1}, e_{R,2}, e_{R,4}, e_{R,11}, e_{R,12}, e_{R,14}, e_{L,2}, e_{P,1}, e_{A,1}\}$$

The transition functions are

$$f_{S,2}(q_{S,2,1}, e_{C,2,1}) = q_{S,2,2},\ f_{S,2}(q_{S,2,2}, e_{R,2}) = q_{S,2,3},$$

$$f_{S,2}(q_{S,2,3}, e_{R,11}) = q_{S,2,4},\ f_{S,2}(q_{S,2,4}, e_{L,2}) = q_{S,2,5},$$

$$f_{S,2}(q_{S,2,5}, e_{R,4}) = q_{S,2,6},\ f_{S,2}(q_{S,2,6}, e_{R,12}) = q_{S,2,7},$$

$$f_{S,2}(q_{S,2,6}, e_{R,14}) = q_{S,2,8},\ f_{S,2}(q_{S,2,7}, e_{A,1}) = q_{S,2,1},$$

$$f_{S,2}(q_{S,2,8}, e_{C,3,1}) = q_{S,2,9},\ f_{S,2}(q_{S,2,9}, e_{P,1}) = q_{S,2,10},$$

$$f_{S,2}(q_{S,2,10}, e_{C,3,1}) = q_{S,2,11}, \text{ and } f_{S,2}(q_{S,2,11}, e_{A,1}) = q_{S,2,1}.$$

The active events sets are

$$\mathbb{H}_{S,2}(q_{S,2,1}) = \{e_{C,2,1}\},\ \mathbb{H}_{S,2}(q_{S,2,2}) = \{e_{R,2}\},$$

$$\mathbb{H}_{S,2}(q_{S,2,3}) = \{e_{R,11}\},\ \mathbb{H}_{S,2}(q_{S,2,4}) = \{e_{L,2}\},$$

$$\mathbb{H}_{S,2}(q_{S,2,5}) = \{e_{R,4}\},\ \mathbb{H}_{S,2}(q_{S,2,6}) = \{e_{R,12}, e_{R,14}\},$$

$$\mathbb{H}_{S,2}(q_{S,2,7}) = \{e_{A,1}\},\ \mathbb{H}_{S,2}(q_{S,2,8}) = \{e_{C,3,1}\},$$

$$\mathbb{H}_{S,2}(q_{S,2,9}) = \{e_{P,1}\},\ \mathbb{H}_{S,2}(q_{S,2,10}) = \{e_{C,3,1}\}$$

$$\text{and } \mathbb{H}_{S,2}(q_{S,2,11}) = \{e_{A,1}\}.$$

The initial state is $x_{S,2,0} = q_{S,2,1}$.

The set of marked states is $\mathbb{Q}_{S,2,m} = \mathbb{Q}_{S,2}$.

The language generated by the automaton is

$$\mathbb{L}(\mathbf{S}_2) = \overline{\left(e_{C,2,1}e_{R,2}e_{R,11}e_{L,2}e_{R,4}(e_{R,12} + e_{R,14}e_{C,3,1}e_{P,1}e_{C,3,1})e_{A,1}\right)^*}$$

Finally, the marked language of the automaton is

$$\mathbb{L}_m(\mathbf{S}_2) = \overline{\left(e_{C,2,1}e_{R,2}e_{R,11}e_{L,2}e_{R,4}(e_{R,12} + e_{R,14}e_{C,3,1}e_{P,1}e_{C,3,1})e_{A,1}\right)^*}$$

The supervisor's automaton is shown in Figure 9.

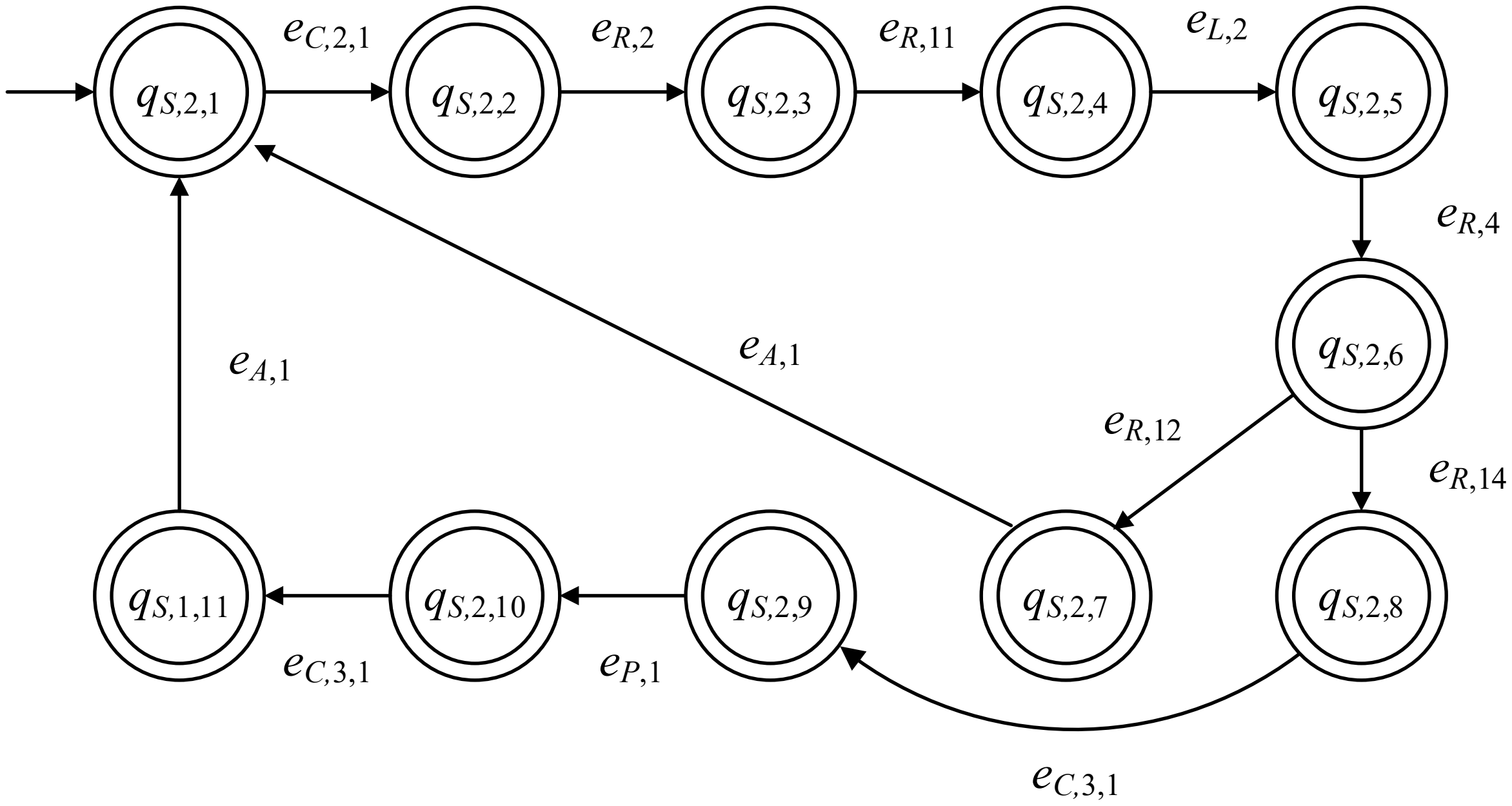

Figure 9: Automaton of $\mathbf{S_2}$

# CHAPTER 5: Conclusions

In this thesis, the mathematical model of a multi-product processing unit was presented. The unit consists of a set of three conveyors, a robot, a lathe, a milling machine, an assembly machine, and a painting machine. Finally, the connection of these elements is made through temporary storage buffers. The products can be divided into two distinct sets depending on the path they will follow through the unit's machines. The mathematical models of the individual subsystems were presented using finite deterministic automata. The two different paths that the products can follow were presented. The desired paths were presented in the form of desired regular languages. The properties of the desired languages with respect to the overall automaton of the system were explored. A supervisory control architecture was designed.